\documentclass[12pt]{amsart}
\usepackage{amssymb}
\usepackage[usenames]{color}

\newtheorem{theorem}{Theorem}

\newtheorem{proposition}{Proposition}
\newtheorem{corollary}{Corollary}

\newtheorem{definition}{Definition}
\newtheorem{lemma}{Lemma}

\begin{document}


\author{V.A.~Vassiliev}
\address{Weizmann Institute of Science, Rehovot, Israel}
 \email{vavassiliev@gmail.com}

\thanks{This work was supported by the Absorption Center in Science of the Ministry of Immigration and Absorption of the State of Israel
}

\title[Equivariant maps of spheres]{Cohomology of spaces of Hopf equivariant maps of spheres}

\begin{abstract}
For any natural numbers $k \leq n$, the rational cohomology ring of the space of continuous maps $S^{2k-1} \to S^{2n-1}$ (respectively, $S^{4k-1} \to S^{4n-1}$), which are equivariant under the Hopf action of the circle (respectively, of the group $S^3$ of unit quaternions),  is naturally isomorphic to that of the Stiefel manifold $V_k({\mathbb C}^n)$ (respectively, $V_k({\mathbb H}^n)$). The natural maps of integral cohomology groups of these spaces of equivariant maps to cohomology of Stiefel manifolds are surjective but not injective.
\end{abstract}

\keywords{Equivariant map, configuration space, Stiefel manifold, spectral sequence}

\subjclass[2020]{57R91, 55R91, 55P91}

\maketitle

\section{Main theorem}

The group $S^1$ (respectively, $S^3$) of complex numbers (respectively, quaternions) of unit norm acts by (left) multiplications on the unit spheres in the complex (respectively, quaternionic) vector spaces. 

\begin{definition} \rm
Denote by $E_{S^1}(k,n)$ the space of continuous maps $f: S^{2k-1} \to S^{2n-1}$ that are equivariant under these actions of $S^1$ on both spheres, i.e. of maps such that $f(\lambda x) =\lambda f(x)$ for any $\lambda \in S^1 \subset {\mathbb C}^1$, $x \in S^{2k-1} \subset {\mathbb C}^k$. 
Denote by $E_{S^3}(k,n)$ the analogous space of continuous maps $f: S^{4k-1} \to S^{4n-1}$ that are equivariant under the left multiplications by quaternions of unit norm. 
\end{definition}

The Stiefel manifold $V_k({\mathbb C}^n)$ of orthonormal $k$-frames in ${\mathbb C}^n$ can be identified with the space of isometric complex linear maps of Hermitian spaces ${\mathbb C}^k \to {\mathbb C}^n$, and, therefore, be considered as a subset of $E_{S^1}(k,n)$. Similarly, the quaternionic Stiefel manifold $V_k({\mathbb H}^n)$ can be considered as a subset of the space $E_{S^3}(k,n)$. 

\begin{theorem}
\label{mainth}
For any natural $k \leq n$, the maps \begin{equation} \label{mt1}H^*(E_{S^1}(k,n), {\mathbb Q})\to  H^*(V_k({\mathbb C}^n), {\mathbb Q}) \ , \end{equation}   \begin{equation} \label{mtm2} H^*(E_{S^3}(k,n), {\mathbb Q}) \to  H^*(V_k({\mathbb H}^n), {\mathbb Q}) \end{equation} induced by inclusions are isomorphisms.

The analogous maps of integer cohomology rings are surjections.
\end{theorem}

\noindent
{\bf Remarks.}
Analogous statements on the injectivity in integer cohomology are not true, see \S \ref{onint}.
 
\noindent
I believe that our method of calculation will be effective in many similar problems.

\section{Preliminaries}
\label{prel}

Recall (see \cite{borel}) that there are ring isomorphisms 
\begin{equation}
\label{borf2}
H^*(V_k({\mathbb C}^n),{\mathbb Z}) \simeq H^*(S^{2n-1} \times S^{2(n-1)-1} \times \dots \times S^{2(n-k+1)-1}, {\mathbb Z}) \ ,
\end{equation}
\begin{equation}
\label{by}
H^*(V_k({\mathbb H}^n),{\mathbb Z}) \simeq H^*(S^{4n-1} \times S^{4(n-1)-1} \times \dots \times S^{4(n-k+1)-1}, {\mathbb Z}).
\end{equation}

\begin{definition} \rm A {\em singular $m$-dimensional simplex} in $E_{S^1}(k,n)$ is a continuous map 
\begin{equation}
F: S^{2k-1} \times \Delta^m \to S^{2n-1},
\label{eqchain}
\end{equation}
 where $\Delta^m$ is the standard $m$-dimensional simplex, and all maps $F(\cdot, \varkappa): S^{2k-1} \to S^{2n-1},$ $\varkappa \in \Delta^m$, are of class $E_{S^1}(k,n)$. Singular $m$-simplices in $E_{S^3}(k,n)$ are defined analogously. The left groups in (\ref{mt1}), (\ref{mtm2}) are the homology groups of arising complexes.

Denote by $\tilde E_{S^1}(k,n)$  and   $\tilde E_{S^3}(k,n)$,  respectively, the subspaces in $E_{S^1}(k,n)$ and $E_{S^3}(k,n)$ consisting of {\em $C^\infty$-smooth} equivariant maps. Singular simplices in $\tilde E_{S^1}(k, n)$ are defined as the maps $($\ref{eqchain}$)$ that are $C^\infty$-smooth {\em in both  variables}, and analogously for   $\tilde E_{S^3}(k, n)$. 
\end{definition}

\begin{proposition}[cf. \cite{JTA}, \S 7]
There are natural isomorphisms  \begin{equation} H^*( E_{S^1}(k,n), {\mathbb Z})  \simeq H^*(\tilde E_{S^1}(k,n), {\mathbb Z}), \label{w1}
\end{equation}
\begin{equation} H^*( E_{S^3}(k,n), {\mathbb Z})  \simeq H^*(\tilde E_{S^3}(k,n), {\mathbb Z}) .
\label{w2}
\end{equation}
\end{proposition}

\noindent
{\it Proof.}
By barycentric subdivisions, the homology class of any $m$-dimensional singular cycle in  $E_{S^1}(k,n)$ can be realized by the class of a continuous map of a closed $m$-dimensional simplicial polyhedron $A $ to  $E_{S^1}(k,n)$, i.e. by a map $\gamma: S^{2k-1} \times A \to S^{2n-1}$; we can also assume that the polyhedron $A$ is embedded in some space ${\mathbb R}^b$. By the Weierstrass theorem, we can approximate this map arbitrarily well  by the restriction of a polynomial map $ \Phi: {\mathbb R}^{2k} \times {\mathbb R}^b \to {\mathbb R}^{2n} \setminus 0$ to the compact subset $S^{2k-1} \times A \subset {\mathbb R}^{2k} \times R^b$. The resulting maps $\varphi_\varkappa \equiv \Phi(\cdot, \varkappa): S^{2k-1} \to {\mathbb R}^{2n}$, $\varkappa \in A$, may not be $S^1$-equivariant, but averaging them over the Haar measure, i.e., replacing any such map by the map taking any point $x \in S^{2k-1}$ to 
$$\int_{0}^{2\pi} e^{i\tau} \varphi_\varkappa(e^{-i\tau} x) d\tau \ \ ,$$
we obtain a map $S^{2k-1} \times A \to {\mathbb R}^{2n} \setminus 0$ equivariant with respect to the action of $S^1$ on $S^{2k-1}$ and ${\mathbb R}^{2n} = {\mathbb C}^n$. This map still approximates well the initial map $\gamma$. Composing it with the standard projection ${\mathbb R}^{2n} \setminus 0 \to S^{2n-1}$, we obtain  a cycle in  $E_{S^1}(k,n)$
 homologous to initial one and representing an element of the right-hand group of (\ref{w1}). In the same way, any homology in $E_{S^1}(k,n)$ between two such cycles in $\tilde E_{S^1}(k,n)$ can be approximated by a homology within $\tilde E_{S^1}(k,n)$, which proves the equality of the classes of these cycles in this right-hand part 
$H^*(\tilde E_{S^1}(k,n), {\mathbb Z})$. The proof of isomorphism (\ref{w2}) is analogous. \hfill $\Box$ 

\begin{proposition}
All maps of class $\tilde E_{S^1}(k,n)$  $($respectively,  $\tilde E_{S^3}(k,n))$ are homotopic to each other within this class.
\end{proposition}

\noindent
{\it Proof.} The space of all equivariant maps $S^{2k-1} \to {\mathbb C}^n \setminus 0$ is obviously homeomorphic to the direct product  $\tilde E_{S^1}(k, n) \times \mbox{Map}({\mathbb C}P^{k-1} \to {\mathbb R}_{> 0})$, in particular it is homotopy equivalent to its subset (and deformation retract) $\tilde E_{S^1}(k, n)$.
 Any two maps $f$, $\tilde f$ of the last space  can be connected in the space of equivariant maps $S^{2k-1} \to {\mathbb C}^n $ by the segment consisting of maps $f_t, t \in [0,1]$, sending any point  $x \in S^{2k-1}$ to $t f(x) + (1-t) \tilde f (x)$. This segment does not lie in the space of maps to $  {\mathbb C}^n \setminus 0$ only if there exists at least one point $x \in S^{2k-1}$ such that the segment $[f(x), \tilde f(x)]$ contains the origin of ${\mathbb C}^n$. This situation can always be removed by an arbitrarily small (and hence keeping the homotopy class) perturbation of the map $\tilde f$ in the class of equivariant maps. Indeed, let us choose an affine chart in ${\mathbb C}P^{k-1}$; let $T \subset S^{2k-1}$ be the image of a cross-section of the Hopf bundle of $S^{2k-1}$ over this chart. For any point $x \in T$ the condition of collinearity of vectors $f(x) $ and $\tilde f(x)$ has codimension $2n-1$, therefore we can arbitrarily slightly perturb our function $\tilde f|_T$ in $2k-2$ variables, keeping it outside a very large ball in $T$, in such a way that this condition does not hold in interior points of this ball. Then this perturbation uniquely extends by equivariantness condition to entire union of Hopf circles over our chart. The new equivariant map is obviously homotopic to $\tilde f$, but can have bad points only over a thin neighborhood of the subspace ${\mathbb C}P^{k-2}$. Choosing another affine chart and repeating the procedure, we get rid of such points outside the preimage of a thin neighborhood of ${\mathbb C}P^{k-3}$. After the $k$-th step we remove all such points and get a map definitely homotopic to both $f$ and $\tilde f$.

The proof for $\tilde E_{S^3}(k,n)$ is completely analogous.\hfill $\Box$ \medskip

\noindent
{\bf Remark.} By a minor additional effort one can prove that the space $\tilde E_{S^1}(k, n)$ is not only connected, but also $2(n-k)$-connected.

\section{Surjectivity of the map (\ref{mt1})}
\label{surje}

Let us fix an orthonormal frame $\{e_1, \dots, e_k\}$ in the Hermitian space ${\mathbb C}^k$ and, consequently, a complete flag ${\mathbb C}^1 \subset \dots \subset {\mathbb C}^{k-1} \subset {\mathbb C}^k$, where any subspace ${\mathbb C}^j$ is spanned by vectors $e_1, \dots, e_j$. Consider ${\mathbb C}^k$ as a subspace in ${\mathbb C}^n$, so that the source $S^{2k-1} \subset {\mathbb C}^k$ of our equivariant maps is a subset of their target space $S^{2n-1} \subset {\mathbb C}^n$. For any $j=1, 2, \dots, k,$ denote by $\Xi_j$ the subset in $\tilde E_{S^1}(k,n)$ consisting of all $S^1$-equivariant smooth maps $S^{2k-1} \to S^{2n-1}$ which are identical on some {\it Hopf circle} in $S^{2k-1} \cap {\mathbb C}^j$ (i.e. on a circle cut from $S^{2k-1}$ by a complex line through the origin in ${\mathbb C}^j \subset {\mathbb C}^k$). 

We will prove that intersection indices with subsets $\Xi_j$ define some cohomology classes of the space $\tilde E_{S^1}(k, n)$, and the restrictions of these classes $[\Xi_i]$ to the subset $V_k({\mathbb C}^n)$ generate multiplicatively the ring $H^*(V_k({\mathbb C}^n),{\mathbb Z})$.

\begin{definition}
\label{defreg}
\rm 
 A map $f \in \tilde E_{S^1}(k, n)$ is {\em regular at a Hopf circle in} $S^{2k-1}$ on which $f$ is identical if for any  point ${\bf x}$ of this circle and any $(2k-2)$-dimensional transversal slice of this circle in $S^{2k-1}$ at this point ${\bf x}$, the differential of the map \ $f - \mbox{\rm Id}$ \ from this slice to ${\mathbb C}^n$ is injective (i.e. has rank $2k-2$) at ${\bf x}$. A map $f \in \Xi_j$ is {\em generic in $\Xi_j$} if it is regular at
any Hopf circle in $S^{2k-1} \cap {\mathbb C}^j$, on which $f$ is identical. (A map $f$ can be generic in $\Xi_j$ and not generic in some $\Xi_l, l>j$.)
\end{definition}

\begin{lemma} 
If a map $f$ is generic in $\Xi_j$ then it is identical on only finitely many Hopf circles in $S^{2k-1} \cap {\mathbb C}^j$.
\end{lemma}

\noindent
{\it Proof.} By definition, Hopf circles on which $f$ is identical are isolated, and our assertion follows from the compactness of the space ${\mathbb C}P^{j-1}$ of all Hopf circles. \hfill $\Box$ \medskip

Define the {\it equivariant $1$-jet space $J^1_{S^1}(S^{2k-1}, S^{2n-1})$} as
the quotient space of the space of $1$-jets of {\em  $S^1$-equivariant} maps $S^{2k-1} \to S^{2n-1}$ by the natural action of $S^1$. In local terms, let us
arbitrarily cover the space ${\mathbb C}P^{k-1}$ by open domains, over any of which the Hopf fibration is trivializable, and choose a smooth cross-section of this fibration over any such domain, i.e. a transversal slice of all Hopf circles over it. Then the open charts of the manifold $J^1_{S^1}(S^{2k-1}, S^{2n-1})$ are spaces of $1$-jets of arbitrary smooth maps of these slices to $S^{2n-1}$. The $1$-jet of such a map at some point of the slice determines uniquely the $1$-jets of any $S^1$-equivariant map $S^{2k-1} \to S^{2n-1}$ defining this jet at all points of the Hopf circle containing this point. In particular, it determines similar $1$-jets of induced maps of all similar slices associated with other  domains of our cover of  ${\mathbb C}P^{k-1}$ at the intersection points of these slices with our Hopf circle; this defines the gluing of our charts of $J^1_{S^1}(S^{2k-1}, S^{2n-1})$ associated with different domains and corresponding slices. The {\it equivariant $0$-jet space $J^0_{S^1}(S^{2k-1}, S^{2n-1})$} is defined in exactly the same way; it is equal to the usual quotient space of $S^{2k-1} \times S^{2n-1}$ by the simultaneous action of the circle on the factors.

Denote by $\Xi_j^{(1)}$ the subset in $J^1_{S^1}(S^{2k-1}, S^{2n-1})$ consisting of equivariant $1$-jets of maps $f: S^{2k-1} \to S^{2n-1}$ such that $f(x) = x$ for all points $x \in S^{2k-1} \cap {\mathbb C}^j$ at which this jet is attached. In particular, any set $\Xi_j^{(1)}$ is the preimage of 
the analogously defined subset $\Xi_j^{(0)} \subset J^0_{S^1}(S^{2k-1}, S^{2n-1})$ under the obvious projection $J^1_{S^1}(S^{2k-1}, S^{2n-1}) \to J^0_{S^1}(S^{2k-1}, S^{2n-1})$. 
Denote by $\Sigma (\Xi_j^{(1)}) $ the subset in $\Xi_j^{(1)}$ consisting of $1$-jets of maps that are not regular (in the sense of Definition \ref{defreg}) at the Hopf circles at which these jets are attached. 

\begin{lemma}
\label{mnlem}
1. $\Xi_j^{(1)}$ is a smooth algebraic submanifold of codimension $(2n-1)+2(k-j)$ in $J^1_{S^1}(S^{2k-1}, S^{2n-1})$.

2. The normal bundle of $\Xi_j^{(1)}$  in $J^1_{S^1}(S^{2k-1}, S^{2n-1})$ has a canonical orientation.

3. $\Sigma (\Xi_j^{(1)}) $ is an algebraic subset in $J^1_{S^1}(S^{2k-1}, S^{2n-1})$ of codimension at   $(2n-1) +2+2(n-j)$.
\end{lemma}

\noindent
{\it Proof.} Statements 1 and 3 follow directly from the definitions. The orientation assumed in statement 2 is induced by a fixed orientation of $S^{2n-1}$ and the complex orientations of the spaces ${\mathbb C}P^{k-1}$ and ${\mathbb C}P^{j-1}$.
\hfill $\Box$

\begin{corollary}
\label{mlem}
For any $j=1, \dots, k$,
intersection indices with subset $\Xi_j \subset \tilde E_{S^1}(k,n)$
 define an element $[\Xi_j]$ of the group $H^{2n-2j+1}(\tilde E_{S^1}(k,n), {\mathbb Z})$. 
\end{corollary}

\noindent
{\it Proof.} For any $m$ consider  the subsets $\Xi_j^{(1)}[m] \supset \Sigma (\Xi_j^{(1)})[m]$ of the space of $1$-jets of $S^1$-equivariant maps (\ref{eqchain}) consisting of $1$-jets, whose restrictions to the subspaces $S^{2k-1} \times \varkappa$, $\varkappa \in \Delta^m$, belong respectively to the subsets $\Xi_j^{(1)}$ and $ \Sigma (\Xi_j^{(1)})$. 
Obviously, all the assertions of Lemma \ref{mnlem} remain valid if we replace the spaces $\Xi_j^{(1)}$ and $ \Sigma (\Xi_j^{(1)})$ respectively by $\Xi_j^{(1)}[m]$ and $ \Sigma (\Xi_j^{(1)})[m]$ and the jet space $J^1_{S^1}(S^{2k-1}, S^{2n-1})$ by the space of all such equivariant $1$-jets of maps (\ref{eqchain}). 

Calculating the homology groups of $\tilde E_{S^1}(k,n)$, we can
replace the complex generated by all smooth singular simplices (\ref{eqchain}) with its  subcomplex generated by such simplices  whose $1$-jet extensions are transversal to the pair of varieties $(\Xi_j^{(1)}[m] , \Sigma (\Xi_j^{(1)})[m] )$. Indeed, by Thom transversality theorem  (the proof \cite{GG} of which is local and therefore works in the equivariant situation), any cycle can be slightly perturbed so that the jet extensions of its restrictions on all simplices become transversal to these pairs as maps of stratified varieties (see \cite{GM}, \S I.1.3).  Such perturbations exist also for the chains defining homology between transversal cycles. 

 For any smooth $(2n-2j+1)$-dimensional singular simplex $($\ref{eqchain}$)$ satisfying this transversality condition, the preimage of the set $\Xi_j^{(1)}$ in $S^{2k-1} \times \Delta^{2n-2j+1}$ consists of finitely many circles of the form $C_\alpha \times \varkappa_\alpha$, where $C_\alpha$ are some Hopf circles in $S^{2k-1} \cap {\mathbb C}^j$ and $\varkappa_\alpha$ are interior points of $\Delta^{2n-2j+1}$, supplied with canonical coorientations. Counting these circles with signs defined by the comparison of these coorientations with the canonical orientation of ${\mathbb C}^{j-1} \times \Delta^{2n-2j+1}$, we obtain the value of our cochain $[\Xi_j]$ on this singular simplex. Analogously, the $1$-jet extension of any smooth singular simplex $($\ref{eqchain}$)$ with $m= (2n-2j+2)$ satisfying this transversality condition does not meet the space $\Sigma (\Xi_j^{(1)})[2n-2j+2]$ of non-regular jets. Therefore, the preimage of the space $\Xi_j^{(1)}[2n-2j+2]$ is a finite union of one-parametric families of Hopf circles in $S^{2k-1} \times \Delta^{2n-2j+2}$, parametrized by one-dimensional closed regular submanifolds in ${\mathbb C}P^{j-1} \times \Delta^{2n-2j+2}$. Moreover, by transversality condition the endpoints of these 1-manifolds in ${\mathbb C}P^{j-1} \times \partial \Delta^{2n-2j+2}$ lie over the maximal faces of $\partial \Delta^{2n-2j+2}$. The number of these endpoints counted with the signs defined by standard coorientations is equal to 0, which proves that our cochain $[\Xi_j]$ is a cocycle. \hfill $\Box$ 

\begin{proposition}
\label{mlem2}
The restrictions of the $k$ cohomology classes \\
$[\Xi_j] \in H^{2(n-j)+1}(\tilde E_{S^1}(k,n), {\mathbb Z})$, $j=1, \dots, k$, to the subspace $V_k({\mathbb C}^n) \subset \tilde E_{S^1}(k,n)$ multiplicatively generate the ring $H^*(V_k({\mathbb C}^n), {\mathbb Z})$.
\end{proposition}

\noindent
{\it Proof.} For any $j=0, 1, \dots, k$ define the subset $A_{j} \subset V_k({\mathbb C}^n)$ as the set of isometric ${\mathbb C}$-linear maps ${\mathbb C}^k \to {\mathbb C}^n$ sending the first $k-j$ vectors $e_1, \dots, e_{k-j}$ of the basic frame to $i e_1$, \dots, $i e_{k-j}$ respectively. Obviously, any $A_{j}$ is diffeomorphic to $V_j({\mathbb C}^{n-k+j}),$ and we have the inclusions $\{i \cdot \mbox{\rm Id}({\mathbb C}^k)\} \equiv A_0 \subset A_1   \subset A_2 \subset \dots \subset A_{k-1} \subset A_k \equiv V_k({\mathbb C}^n), $ where any $A_{j}$ is embedded into $A_{j+1}$ as a fiber of the fiber bundle $A_{j+1} \to S^{2(n-k+j)+1}$ taking any map ${\mathbb C}^k \to {\mathbb C}^n$  into the image of the basic vector $e_{k-j}$ under this map. The manifolds $ A_{j+1}$ and $\Xi_{k-j} $ meet transversally in $\tilde E_{S^1}(k,n)$, and their intersection is equal to a different fiber of the fiber bundle  $A_{j+1} \to S^{2(n-k+j)+1}$, namely the fiber over the vector $e_{k-j}$ itself. Indeed, any point of this intersection is a linear map of class $V_k({\mathbb C}^n)$ which, on one hand, multiplies all vectors $e_1, \dots, e_{k-j-1}$ by $i$, and on the other 
sends some vector from the linear span of vectors $e_1, \dots, e_{k-j}$ to itself.
Therefore the restriction of the cohomology class $[\Xi_{k-j}]$ to $A_{j}$ (respectively, to $A_{j+1}$) is equal to zero (respectively, coincides with the class induced from the basic cohomology class of $S^{2(n-k+j)+1}$ by the projection of this fiber bundle). Proposition \ref{mlem2} follows from this by induction over the spectral sequences of these fiber bundles and the multiplicativity of these spectral sequences, see for example, \cite{FF}, \S 24.3. \hfill $\Box$ 

\begin{corollary}
\label{mcor}
The map 
$$H^*(E_{S^1}(k,n), {\mathbb Q})\to  H^*(V_k({\mathbb C}^n), {\mathbb Q}) $$ defined by the identical embedding is epimorphic. \hfill $\Box$
\end{corollary}

\section{Injectivity of map (\ref{mt1})}

For any topological space $X$ and natural number $t$, denote by $B(X, t)$ the $t$-configuration space of $X$ (that is, the space of $t$-element subsets of $X$ with a natural topology). Denote by $\pm {\mathbb Z}$ the {\em sign local system} of groups on $B(X, t)$: it is locally isomorphic to ${\mathbb Z}$, but loops in $B(X, t)$ act on its fibers as multiplication by $1$ or $-1$ depending on the parity of permutations of $t$ points defined by these loops. $\bar H_*$ denotes the Borel--Moore homology group (that is, the homology group of the complex of locally finite singular chains).

\begin{proposition}
\label{dva}
There is a spectral sequence $E_r^{p,q}$ converging to the group $\tilde H^*(E_{S^1}(k,n),{\mathbb Z})$, whose term $E_1^{p,q}$ is trivial if $p\geq 0$ and is defined by the formula
\begin{equation}
\label{mfo}
  E_1^{p,q} \simeq       \bar H_{-2pn-q}(B({\mathbb C}P^{k-1},-p), \pm {\mathbb Z})
\end{equation}
for $p<0$.
\end{proposition}

\noindent {\it Proof.} This is a direct application of Theorem 3 of \cite{JTA}. Namely, the elements  $X$, ${\mathbb R}^a$, $G$, $S^{W-1}$ and ${\mathbb R}^W$ of this theorem appear in our case  respectively as  $S^{2k-1}$, ${\mathbb C}^k$, $S^1$, $S^{2n-1},$ and ${\mathbb C}^n$; the subset $\Lambda \subset S^{W-1}$ is empty,  and $C\Lambda$ is the origin in ${\mathbb C}^n$. \hfill $\Box$

\begin{proposition}[see \cite{how}, Lemma 2(B)] 
\label{tri}
 For any $j \geq 0$ and $t \geq 1$ there is a group isomorphism
\begin{equation}
\label{cpnk} \bar H_j(B({\mathbb C}P^{k-1}, t), \pm {\mathbb Z} \otimes {\mathbb Q}) \simeq H_{j-t(t-1)}(G_t({\mathbb C}^{k}),{\mathbb Q}),
\end{equation}
where
$G_t({\mathbb C}^{k})$ is the Grassmann manifold of $t$-dimensional complex subspaces in
${\mathbb C}^{k}$.
In particular, the left-hand group $($\ref{cpnk}$)$ is trivial if $t>k$. \hfill $\Box$
\end{proposition} 

By the propositions \ref{dva} and \ref{tri}, the total dimension of the group $\tilde H^*(\tilde E_{S^1}(k,n), {\mathbb Q})$ does not exceed 
$$\sum_{t=1}^k  \dim H^*(G_t({\mathbb C}^k), {\mathbb Q}) = \sum_{t=1}^k\binom{k}{t} = 2^k-1, $$ 
which by (\ref{borf2}) is equal to the total dimension of the reduced homology group of the Stiefel manifold  $V_k({\mathbb C}^n)$. So, the map (\ref{mt1}) is an epimorphic (by Corollary \ref{mcor}) map, the dimension of whose source does not exceed that of the target; hence it is an isomorphism. \hfill $\Box$

\section{Bijectivity of the map (\ref{mtm2})}

The assertions of Theorem \ref{mainth} on the map (\ref{mtm2}) and its analog over the integer coefficients can be proved in exactly the same way. In particular, the first term of the spectral sequence calculating the group $\tilde H^*(\tilde E_{S^3}(k,n),{\mathbb Z})$ and analogous to (\ref{mfo}) is given by 
\begin{equation} {\bf E}_1^{p, q} \simeq  \bar H_{-4pn-q}(B({\mathbb H} P^{k-1}, -p), \pm {\mathbb Z}) 
\label{mfoq}
\end{equation}
for $p< 0$ and is trivial for $p \geq 0$; 
 the quaternionic version of equality (\ref{cpnk}) is 
\begin{equation}
\bar H_j(B({\mathbb H} P^{k-1}, t), \pm {\mathbb Z} \otimes {\mathbb Q}) \simeq H_{j-2t(t-1)}(G_t({\mathbb H}^{k}), {\mathbb Q}) \ .
\end{equation}

\section{The first column of spectral sequences}

By formula (\ref{mfo}), the groups $E_1^{-1, q}$ of our spectral sequence are equal to ${\mathbb Z}$ for $q = 2n-2k+ 2, 2n-2k +4, \dots, 2n$, and are trivial for all other $q$.

\begin{proposition}
The integral cohomology classes  of dimensions $2n-2j+1$, $j = 1, \dots, k$, of the space $\tilde E_{S^1}(k,n)$, which  correspond to the generators of the groups $E_1^{-1, 2n-2j+2}$,  coincide with the classes $[\Xi_j]$ considered in \S \ref{surje}.
\end{proposition}

\noindent
{\it Proof.} By the construction of \cite{JTA}, the basic element of  the group
$E_1^{-1, 2n - 2j +2} $ defines the $(2n-2j+1)$-dimensional cohomology class of the space of equivariant maps $S^{2k-1} \to {\mathbb C}^n \setminus \{0\}$ equal to the linking number in the vector space of all  equivariant maps $S^{2k-1} \to {\mathbb C}^n$ with the subset consisting of  maps which send to 0 some Hopf circle from the sphere $S^{2k-1} \cap {\mathbb C}^j$. (See \S 7 of \cite{JTA} for the technique of finite-dimensional approximations justifying the argumentation in terms of linking and intersection numbers in spaces of equivariant maps.) This subset of singular maps is the boundary of the set of maps acting on some such circle not as the multiplication by 0 but as the multiplication by an arbitrary positive number. The intersection of this set with the subspace of maps with target $S^{2n-1} \subset {\mathbb C}^n \setminus 0 $ is exactly the variety $\Xi_j$. \hfill $\Box$

\begin{corollary}
All maps $E_1^{-1, q} \to E_{\infty}^{-1, q}$ of our spectral sequences are isomorphisms. \hfill $\Box$
\end{corollary}

\section{On the integer homology}
\label{onint}

The following example proves that the assertions of Theorem \ref{mainth} on the injectivity of the maps (\ref{mt1}), (\ref{mtm2}) cannot be extended to  cohomology with integer coefficients, in particular they are not reflections of homotopy equivalences.
\begin{proposition}
\label{cont}
1. The group 
$H^{4n-4k+2}(\tilde E_{S^1}(k,n), {\mathbb Z})$ 
is trivial, but the group \\
$H^{4n-4k+2}(\tilde E_{S^1}(k,n), {\mathbb Z}_2)$ is isomorphic to ${\mathbb Z}_2$. 

2. The group $H^{8n-8k+6}(\tilde E_{S^3}(k,n), {\mathbb Z}) $ 
is trivial, but $H^{8n-8k+6}(\tilde E_{S^3}(k,n), {\mathbb Z}_2)$ is isomorphic to ${\mathbb Z}_2$.
\end{proposition}

\begin{corollary}
The groups $H^{4n-4k+3}(\tilde E_{S^1}(k,n), {\mathbb Z})$ and $H^{8n-8k+7}(\tilde E_{S^3}(k,n), {\mathbb Z}) $ contain a non-trivial 2-torsion. \hfill $\Box$
\end{corollary}

\begin{lemma}
For any natural $k \leq n$ there exists a spectral sequence  $e^{p,q}_r    $ calculating the mod 2 cohomology of the space $E_{S^1}(k,n)$; its groups $e_1^{p,q}$  are trivial for $p \geq 0$ and
\begin{equation}
\label{mfo2}
e_1^{p,q} \simeq       \bar H_{-2pn-q}(B({\mathbb C}P^{k-1},-p),  {\mathbb Z}_2)
\end{equation}
for $p<0$.
\end{lemma}

\noindent
{\it Proof.} The proof is a literal reduction modulo 2 of the proof of Proposition \ref{dva} given in \cite{JTA}. \hfill $\Box$

\begin{lemma} 
\label{le47}
For any $k \leq n$, 

1. The groups $E_1^{-2, q} $ 
are trivial for all $q \leq 4(n-k)+4$. The groups $e_1^{-2,q}$ are trivial for all 
$q < 4(n-k)+4   $ and equal to ${\mathbb Z}_2$ for $q = 4(n-k)+4$.

2. The minimal total dimension $p+q$ of non-trivial groups $E_1^{p,q}$ and $e_1^{p,q}$ with $p \leq -3$ is at least $3(2n-2k+1)$.  
\end{lemma}

\noindent
{\it Proof.} 1. By (\ref{mfo2}),  $e_1^{-2, 4(n-k+1)} = \bar H_{4(k-1)}(B({\mathbb C}P^{k-1}, 2),  {\mathbb Z}_2)$. $B({\mathbb C}P^{k-1}, 2)$ is a connected $4(k-1)$-dimensional manifold, hence its $4(k-1)$-dimensional Borel--Moore mod 2 homology group  is isomorphic to ${\mathbb Z}_2$. This manifold is not $\pm {\mathbb Z}$-orientable since any loop in it permuting two points changes this orientation. Therefore, the group
$ E_1^{-2, 4(n-k+1)} \simeq \bar H_{4(k-1)}(B({\mathbb C}P^{k-1}, 2), \pm {\mathbb Z})$ is trivial. All  homology groups of higher dimensions of this manifold (corresponding to groups $E_1^{-2,q} and $ $e_1^{-2, q}$ with smaller $q$) are obviously trivial. 

Statement 2 follows immediately from the formulas (\ref{mfo}) and (\ref{mfo2}) and from the fact that $B({\mathbb C}P^{k-1}, -p)$ is a $-2p(k-1)$-dimensional manifold. 
\hfill $\Box$
\medskip

\noindent 
{\it Proof of Proposition \ref{cont}.} The differential $d^1: e_1^{-2, 4(n - k +1)} \to e_1^{-1, 4(n - k +1)} $ is trivial because  the last group is (depending on $n$ and $k$) either trivial or generated by the reduction mod 2 of a cycle defining an integer cohomology class taking arbitrary integer values on cycles of dimension $4(n-k)+3$. By statement 2 of Lemma \ref{le47}, all higher differentials that hit the group $e_r^{-2,4(n-k+1)} \simeq {\mathbb Z}_2$ are also trivial. Also by Lemma \ref{le47}, all groups $E_1^{p,q} $ with $p+q = 4n-4k+2$ are trivial. This implies  statement 1  of Proposition \ref{cont}. The proof of statement 2 is completely analogous.
\hfill $\Box$ \medskip

\noindent
{\bf Remark.}
The non-trivial element of the group $H^{4n-4k+2}(\tilde E_{S^1}(k,n), {\mathbb Z}_2)$  has the meaning of the intersection index with the set of maps which are identical on two different Hopf circles. The subvariety $V_k({\mathbb C}^n)  \subset  \tilde E_{S^1}(k,n)$ is in a very non-general position with this set.

\end{document}